\documentclass[12pt,leqno]{article}
\usepackage{amsmath}
\usepackage{amssymb}
\newtheorem{theorem}[equation]{Theorem}
\newtheorem{proposition}[equation]{Proposition}
\newtheorem{corollary}[equation]{Corollary}
\newtheorem{lemma}[equation]{Lemma}

\newtheorem{definition}[equation]{{\it Definition}\rm }
\newtheorem{remark}[equation]{\it Remark}
\newcommand{\NO}{NO$\frac{84}{90}$}
\newcommand{\rankcDi}{\mathrm{rank}_{\mathbf{Z}}\{c_1(D_i)\}_{i=1}^l}
\numberwithin{equation}{section}

\title{Holomorphic Curves and Integral Points off Divisors\\
}
\author{Junjiro Noguchi and J\"org Winkelmann}
\begin{document}
\setlength{\baselineskip}{18pt}
\maketitle
\begin{abstract}
We deal with the distributions of holomorphic
curves and integral points off divisors.
We will simultaneouly prove an optimal dimension estimate
from above of a subvariety $W$ off a divisor $D$ which
contains a Zariski dense entire holomorphic curve,
or a Zariski dense $D$-integral point set,
provided that in the latter case everything is defined over a number field.
Then, if the number of components of $D$ is large,
the estimate leads to the constancy of such a holomorphic curve
or the finiteness of such an integral point set.
At the begining, we extend logarithmic Bloch-Ochiai's Theorem
to the K\"ahler case.
\end{abstract}

\section{Introduction and results.}

There have been found a number of interesting analogues
among three themes, the Nevanlinna theory, the Kobayashi hyperbolicity,
and the Diophantine approximation.
Cf., e.g., [Ko70], [Ko98], [L87], [L91], [N91], [N92], [N97],
[\NO], [V96].
In this paper we deal with the distributions of holomorphic
curves and integral points off divisors.
Let $V$ be a complex projective algebraic variety and
let $D$ be a hypersurface (an effective reduced divisor)
of $V$.
We will prove an optimal dimension estimate from above
of a subvariety $W \subset V$ off a divisor $D$ which
admits, roughly speaking, an entire holomorphic
curve $f:\mathbf{C} \to W\setminus D$ with Zariski dense image,
or which contains a Zariski dense $D$-integral point set,
provided that everything is defined over a number field.
Then, if the number of components of $D$ is large,
the estimate leads to the constancy of such a holomorphic curve
or the finiteness of such an integral point set.
Here we have the following counter objects:
\begin{align}
\label{1.0}
\hbox{a non-constant}&\hbox{ holomorphic curve } f:\mathbf{C} \to V \setminus D \\
& \Longleftrightarrow
\hbox{ an infinite }D\hbox{-integral point set of }V. \nonumber
\end{align}

As for holomorphic curves, we deal with them in a compact
K\"ahler manifold unless the algebraicity is needed.
For we think it is the most suited place for holomorphic curves.

Let $M$ be a compact complex manifold of dimension $m$ and
let $\{D_i\}_{i=1}^l$ be a family of hypersurfaces of $M$.
We say that $\{D_i\}_{i=1}^l$ is {\it in general position} if for any
distinct indices $1\leqq i_1,\ldots, i_k \leqq l$,
the codimension of every irreducible component of the intersection
$\bigcap_{j=1}^k D_{i_j}$ is $k$ for $k\leqq m$, and
$\bigcap_{j=1}^k D_{i_j}=\emptyset$ for $k>m$.
This notion is defined for singular $M$ as well.
Assume that $M$ is K\"ahler.
We denote by $c_1(D_i)\in H^2(M,\mathbf{Z})_i\cap H^{1,1}(M, \mathbf{R})$
the Poincar\'e dual of $D_i\in H_{2n-2}(M,\mathbf{Z})_i$,
where $H^2(M,\mathbf{Z})_i$ (resp., $H_{2n-2}(M,\mathbf{Z})_i$) stands
for the image of $H^2(M,\mathbf{Z})$ (resp., $H_{2n-2}(M,\mathbf{Z})$)
in the space $H^2(M,\mathbf{R})$ (resp., $H_{2n-2}(M,\mathbf{R})$).
Let $\rankcDi$ denote the $\mathbf{Z}$-rank of
the subgroup of $H^2(M,\mathbf{R})$
generated by $\{c_1(D_i)\}_{i=1}^l$.
For a subvariety $W \subset M$ we denote by $q(W)$ the irregularity
of $W$ ($=$ the dimension of the vector space of global holomorphic
1-forms of a desingularization model of $W$, which is independent of
the model), and by $\#\{W \cap D_i\not=W\}$ the number of distinct
non-empty cut loci $W \cap D_i, 1 \leqq i \leqq l$,
such that $D_i\cap W \not= W$.

In the present paper it is assumed that an ample or a very ample divisor
means an ample or a very ample effective one.
For a divisor $D$ we write the same notation $D$ for its support
unless confusion occurs.
We use the standard notations,
$\mathbf{Z}, \mathbf{Q}, \mathbf{R}$, and $\mathbf{C}$
for the sets of integers, rational numbers, real numbers,
and complex numbers, respectively.

The first purpose of this paper is to prove the following.
\medskip

\begin{theorem}
\label{1.1}
Let $M$ be a compact K\"ahler manifold of dimension $m$.
Let $\{D_i\}_{i=1}^l$ be a family of hypersurfaces of $M$
in general position.
Let $W \subset M$ be a subvariety such that there is a
non-constant holomorphic curve
$f:\mathbf{C} \to W \setminus\underset{D_i \not\supset W}\bigcup D_i$
with Zariski dense image.
Then we have that
\begin{enumerate}
\item
$\#\{W \cap D_i\not=W\} + q(W) \leqq \dim W + \rankcDi$:
\item
Assume that all $D_i$ are ample.
Then we have
$$
(l-m) \dim W \leqq
m \left( \rankcDi-q(W) \right)^+.
$$
\end{enumerate}
\end{theorem}

Here $(\cdot)^+$ stands for
the maximum of $0$ and the number.

Let $\mathrm{NS}(M)$ denote the
Neron-Severi group of $M$;
i.e., $\mathrm{NS}(M)=\mathrm{Pic}(M)/\mathrm{Pic}^0(M)$.
We know that
\begin{equation}
\label{1.2}
\rankcDi \leqq \mathrm{rank}_{\mathbf{Z}}\,
\mathrm{NS}(M).
\end{equation}

The following corollary provides also examples.
\begin{corollary}
\label{1.3}
{\it
Let the notation be as above.
\begin{enumerate}
\item
Assume that all $D_i$ are ample and that
$l > m(\mathrm{rank}_{\mathbf{Z}}\, \mathrm{NS}(M)+1)$.
Then $M \setminus \bigcup_{i=1}^l D_i$
is complete hyperbolic and hyperbolically imbedded into $M$.
\item
Let $X \subset \mathbf{P}^m(\mathbf{C})$ be an irreducible
subvariety, and let $D_i, 1 \leqq i \leqq l$, be distinct
hypersurface cuts of $X$ that are in general position as
hypersurfaces of $X$.
If $l > 2 \dim X$, then $X \setminus \bigcup_{i=1}^l D_i$
is complete hyperbolic and hyperbolically imbedded into $X$.
\item
Let $\{D_i\}_{i=1}^l$ be a family of ample hypersurfaces
of $M$ in general position.
Let $f:\mathbf{C}\to M$ be a holomorphic curve such that
for every $D_i$, either $f(\mathbf{C})\subset D_i$, or
$f(\mathbf{C})\cap D_i=\emptyset$.
Assume that $l>m$.
Then $f(\mathbf{C})$ is contained in an algebraic subspace $W$ 
of $M$ such that
$$
\dim W \leqq \frac{m}{l-m} \mathrm{rank}_{\mathbf{Z}}\, \mathrm{NS}(M).
$$
In special, if $M=\mathbf{P}^m(\mathbf{C})$, then we have
$$
\dim W \leqq \frac{m}{l-m}.
$$

\end{enumerate}
}
\end{corollary}

The above Corollary \ref{1.3}, (ii) for $X=\mathbf{P}^m(\mathbf{C})$
was given by Babets [B84], but his proof seems to carry some
incompleteness and confusion.
In the case of $\mathbf{P}^m(\mathbf{C})$ and hyperplanes $D_i$,
the above (ii) with $X=\mathbf{P}^m(\mathbf{C})$
and (iii) for
$f:\mathbf{C} \to \mathbf{P}^m(\mathbf{C}) \setminus \bigcup_{i=1}^lD_i$
were first proved by Fujimoto [Fu72] and Green [G72],
where the linearity of $W$
was also proved, and by their examples the dimension estimate
is best possible in general
(see Remark \ref{4.12}).

In the proof of Theorem \ref{1.1}, the logarithmic version of
Bloch-Ochiai's theorem ([Bl26], [O77], [K80]) proved by [N77], [N80],
[N81], [NO$\frac{84}{90}$]) will play a crucial role.
For the proof of Theorem \ref{1.1}, (i),
we will extend logarithmic Bloch-Ochiai's theorem to
the K\"ahler case (cf.\ [N96] and {\bf Theorem \ref{2.11}}).

For a complex torus $M$ these methods
are not effective very much, because $M$
may have a large Neron-Severi group
($\mathrm{rank}_{\mathbf{Z}}\,\mathrm{NS}(M)\leqq (\dim M)^2$).
Here a different method yields better results.

\begin{proposition}
\label{1.4}
Let $M$ be a complex torus and let $m_0$ be the positive minimal
dimension of complex subtori of $M$.
Let $\{D_i\}^l_{i=1}$
be a family of hypersurfaces of $M$ in general position.
Assume that at least one of the following conditions is fulfilled:
\begin{enumerate}
\item
$l >  \dim M - m_0$ and all $D_i$ are ample,
\item
$l > 2(\dim M- m_0)$.
\end{enumerate}
Then
$M \setminus \bigcup_{i=1}^l D_i$ is complete hyperbolic and
hyperbolically imbedded into $M$.
\end{proposition}

As a special case this contains the result of M. Green [G78] that
the complement of a divisor in a simple torus is 
complete hyperbolic and hyperbolically imbedded.
In the case of simple torus $M$, $m_0=\dim M$ and any hypersurface
of $M$ is ample.

In the Diophantine approximation, Theorem \ref{1.1}, (i)
is analogous to Vojta [V96], Corollary 0.3.
Here we apply the analogue of (\ref{1.0}).
Then, by making use of [V96], Theorem 0.2 and Corollary 0.3
in the same way as in the proofs of Theorem \ref{1.1} and
Corollary \ref{1.3},
we have the following results in the Diophantine approximation
(for the notation, see \S3).

\begin{theorem}
\label{1.5}
Assume that everything is  defined over a number field $K$,
and $S$ is a finite subset of a set
$M(K)$ of inequivalent non-trivial places of $K$ with product formula
such that $S$ contains all infinite places.
Let $V$ be a projective smooth variety of dimension $m$.
Let $\{D_i\}_{i=1}^l$ be a family of ample hypersurfaces of
$V$ in general position.
Let $W \subset V$ be a subvariety of $V$.
Assume that there exists a Zariski dense
$(\sum_{D_i \not\supset W} D_i \cap W, S)$-integral point set
of $W(K)$ in $W$.
Then we have
$$
(l-m) \dim W \leqq
m \left( \rankcDi-q(W) \right)^+.
$$
\end{theorem}

\begin{corollary}
\label{1.6}
Let the notation be as above.
\begin{enumerate}
\item
Assume that all $D_i$ are ample and that
$l > m(\mathrm{rank}_{\mathbf{Z}}\, \mathrm{NS}(V)+1)$.
Then any $(\sum_{i=1}^l D_i,S)$-integral point set
of $V(K)$ is finite.
\item
Let $X \subset \mathbf{P}^m_K$ be an irreducible
subvariety, and let $D_i, 1 \leqq i \leqq l$, be distinct
hypersurface cuts of $X$ that are in general position as
hypersurfaces of $X$.
If $l > 2 \dim X$, then any $(\sum_{i=1}^l D_i,S)$-integral point
set of $X(K)$ is finite.
\item
Let $D_i, 1 \leqq i \leqq l$, be ample divisors
of $V$ in general position.
Let $A$ be a subset of $V(K)$ such that
for every $D_i$, either $A \subset D_i$, or
$A$ is a $(\sum_{D_i\not\supset A} D_i,S)$-integral point set.
Assume that $l > m$.
Then $A$ is contained in an algebraic subvariety $W$ of $V$ such that
$$
\dim W \leqq \frac{m}{l-m} \mathrm{rank}_{\mathbf{Z}}\, \mathrm{NS}(V).
$$
In special, if $V=\mathbf{P}^m_K$, then we have
$$
\dim W \leqq \frac{m}{l-m}.
$$
\end{enumerate}
\end{corollary}

The dimension estimates obtained above are optimal
(see Remark \ref{4.12}).
M. Ru and P.-M. Wong [RW91]
dealt with the special case of $V=\mathbf{P}^m_K$ and hyperplanes
$D_i$, and proved that if $A$ is a $(\sum_{i=1}^lD_i,S)$-integral
point set, then $A$ is contained in a finite union $W$ of
linear subspaces such that
$$
\dim W \leqq (2m+1 -l)^+.
$$
Cf.\ $\dim W \leqq m/(l-m)$ of Corollary \ref{1.6}, (iii).

{\it Acknowledgement.}
The first author is grateful to Professor Masa-Hiko Saito for
helpful conversations on Deligne's Hodge theory.
The present joint work was partially supported by
Japan-U.S. Cooperative Science Program, and Japan-U.S.
Mathematical Institute at Johns Hopkins University,
to which the authors express sincere gratitude.
The second author also wishes to thank the ``Schweizer
Nationalfonds'' for partial support.

\section{Logarithmic Bloch-Ochiai's Theorem for K\"ahler\\ manifolds.}

\indent
In this section we assume that
$M$ is a compact K\"ahler manifold of dimension $m$
and let $D$ be a hypersurface of $M$.
We denote by $\Omega_M^1$ (resp., $\Omega_M^1(\log D)$)
the sheaf of (resp., logarithmic) 1-forms (resp., along $D$).
We set
$$
q(M)=\dim H^0(M, \Omega_M^1), \qquad
\bar q(M\setminus D)=\dim H^0(M, \Omega_M^1(\log D))
$$
and call $q(M)$ (resp., $\bar q (M\setminus D)$)
the (resp., {\it logarithmic}) irregularity of $M$
(resp., $M\setminus D$).
Set
$$
M'=M\setminus D.
$$
By the Hodge theory due to Deligne [D71] we have
\begin{equation}
\label{2.1}
\bar q(M')-q(M)  = b_1(M')-b_1(M).
\end{equation}
Furthermore, global logarithmic 1-forms over $M$ are
closed (Deligne [D71], Noguchi [N95]).

We recall the construction of the quasi-Albanese mapping due to
Iitaka [I76].
Let $D_i$ ($i=1,\ldots,d$) denote the connected components of the
smooth part $D\setminus \mathrm{Sing}\, D$ of $D$.
Let $\gamma_i$ denote small circles around $D_i$;
that is, choose a point $p_i\in D_i$ and local coordinates 
$w_k^{(i)}, 1 \leqq k \leqq m$,
about $p_i$ such that $D=\{w_1^{(i)}=0\}$ near $p_i$,
and let
$\gamma_i:S^1=\{e^{i\theta}; \theta \in \mathbf{R} \}\to N$
be given by $e^{i\theta}\mapsto
(\epsilon e^{i\theta},0,\ldots,0)$ with small $\epsilon>0$.
Then the Mayer-Vietoris sequence implies that
$\gamma_i\in H_1(M',\mathbf{Z}), 1 \leqq i \leqq d$,
generate the kernel of the natural map 
$H_1(M',\mathbf{Z})\to H_1(M,\mathbf{Z})$.
If $\omega\in H^0(M,\Omega_M^1(\log D))$ has a pole on $D_i$, then
$\int_{\gamma_i}\omega \not= 0$.
On the other hand, $\int_{\gamma_i}\omega=0$ for all closed 1-forms
$\omega$ which are holomorphic on $D_i$.
Hence equality $(\ref{2.1})$ implies that the
$\mathbf{Q}$-span of $\gamma_i$'s in
$H_1(M',\mathbf{Q})$ is a $\mathbf{Q}$-vector space
of dimension $\rho:=\dim H^0(\Omega_M^1(\log D))
/H^0(\Omega_M^1)$.
Together we obtain the following:
We may renumber $\gamma_i$'s and
choose logarithmic 1-forms $\omega_k$ ($k\in\{1,\ldots,\rho\}$)
with the following properties:

{\samepage
\begin{equation}
\label{2.1.1}
\end{equation}
\begin{enumerate}
\item
$\omega_1,\ldots,\omega_\rho$ form a basis of the vector space
$H^0(M, \Omega^1_M(\log D))/H^0(M, \Omega^1_M)$,
\item 
$\int_{\gamma_j}\omega_k
=\sigma_{jk}\in \mathbf{Z}$ for all $1\leqq i \leqq d$, $1\leqq k\leqq \rho$,
\item
$\sigma_{jk}=\delta_{jk}$ (Kronecker's symbol) if $j\leqq\rho$.
\end{enumerate}}
\noindent
Now, we define the {\it quasi-Albanese variety}
$A(M')$ as the complex quotient Lie group
$H^0(M, \Omega^1_M(\log D))^*/H_1(M',\mathbf{Z})$,
where ``$\empty^*$'' stands for the dual of the vector space,
and the {\it quasi-Albanese mapping}
$$
\alpha: x \in M' \to \left[\omega\mapsto\int_{x_0}^x\omega\right]
\in A(M'),
$$
where $x_0 \in M'$ is a fixed base point.
By (\ref{2.1.1}) the short exact sequence of complex vector spaces
$$
0 \to \left( H^0(M, \Omega^1_M(\log D))/H^0(M, \Omega^1_M)\right)^* 
\to H^0(M, \Omega_M^1(\log D))^*
\to H^0(M,\Omega_M^1)^* \to 0
$$
yields a short exact sequence of complex Lie groups
\begin{equation}
\label{2.2}
1 \to {(\mathbf{C}^*)^\rho}
 \to {A(M')}
\stackrel{\pi}{\to} {A(M)}
\to 1,
\end{equation}
where $\mathbf{C}^*$ denotes the multiplicative group of
non-zero complex numbers, and
$A(M)$ denotes the ordinary Albanese variety (torus) of $M$.

We are now going to compactify $A(M')$. We start with the standard
compactification
$$
\overline{\mathbf{C}^*}=\mathbf{P}^1=\mathbf{C}^*\cup\{0\}\cup\{\infty\},\qquad
(\mathbf{C}^*)^\rho\hookrightarrow
({\overline{\mathbf{C}^*})^\rho}=(\mathbf{P}^1)^\rho.
$$
This is an algebraic compactification.
In addition, it is
equivariant for the multiplicative
$(\mathbf{C}^*)^\rho$-action on itself by translations.
Therefore the $(\overline{\mathbf{C}^*})^\rho$-fiber bundle associated to the 
$(\mathbf{C}^*)^\rho$-principal bundle $A(M')\to A(M)$
via the $(\mathbf{C}^*)^\rho$-action on $(\overline{\mathbf{C}^*})^\rho$ 
yields a compactification $\overline{A(M')}$ of $A(M')$
with projection $\bar\pi: \overline{A(M')} \to A(M)$.
This compactification is obviously $(\mathbf{C}^*)^\rho$-equivariant.
Furthermore, if $\lambda_a$ denotes the automorphism of $A(M)$ given by
a translation with $a\in A(M)$,
then the $(\mathbf{C}^*)^\rho$-principal bundle $\lambda_a^*A \to A(M)$
is isomorphic as $(\mathbf{C}^*)^\rho$-principal bundle to
$A(M')\to A(M)$.
This implies that $\lambda_a$ can be lifted to an automorphism of
$\overline{A(M')}$ which stabilizes $A(M')$ and commutes with
the $(\mathbf{C}^*)^\rho$-principal action on $A(M')$. 
It hence follows that $A(M')\hookrightarrow\overline{A(M')}$
is $A(M')$-equivariant.
\begin{lemma}
\label{2.3}
Let the notation be as above.
Then the quasi-Albanese mapping
$\alpha:M' \to A(M')$ extends to a meromorphic mapping
$\overline\alpha:M \to \overline{A(M')}$ which is
holomorphic on $M\setminus\mathrm{Sing}\, D$.
\end{lemma}

{\it Proof.}
Let $p\in D\setminus\mathrm{Sing}\, D$.
Then there is a $k\in\{1,\ldots,d\}$ such that $p \in D_k$.
Since the bundle $\overline{A(M')} \stackrel{\bar\pi}{\to} A(M)$
is locally trivial,
there is a small open neighborhood $W$ of
$\bar\pi\circ\alpha(p)$ in $A(M)$ such that
$\bar\pi^{-1} W \cong (\mathbf{C}^*)^\rho\times W$.
Let $\beta=\bar\pi\circ\alpha:M \to A(M)$ be the Albanese mapping of $M$.
We may write $\alpha$ near $p$ as
$$
\alpha=(\alpha_1,\ldots,\alpha_\rho, \beta),
$$
where $\alpha_j(x) = \int_{x_0}^x \omega_j$.
The choice of $\omega_j$ now implies that $\lim_{x\to p}\alpha_j(x)$ 
converges in $\mathbf{P}^1$: more precisely, $\lim_{x\to p}\alpha_j(x)=0$
if $\sigma_{jk}>0$, and $\lim_{x\to p}\alpha_j(x)=\infty$ if
$\sigma_{jk}<0$.
If $\sigma_{jk}=0$, then $\omega_j$ is holomorphic on $D_k$
and $\lim_{x\to p}\alpha_j(x)$ has a limit inside $\mathbf{C}^*$.
Hence $\alpha$ extends to a holomorphic mapping from
$M\setminus\mathrm{Sing}\, D \to \overline{A(M')}$.

Now, let $p \in \mathrm{Sing}\, D$.
Observe that $\mathrm{Sing}\, D$ has codimension at least 
two in $M$.
Then there is a small open neighborhood $U$ of $p$ in $M$
such that $U\setminus\mathrm{Sing}\, D$ is simply connected.
The mapping $\alpha:U\setminus \mathrm{Sing}\, D \to \overline{A(M')}$
is lifted to a mapping $\tilde\alpha$ from $U \setminus\mathrm{Sing}\, D$
to the universal covering of $\overline{A(M')}$ which is
$(\mathbf{P}^1)^\rho\times\mathbf{C}^{q(M)}$.
Hence, $\tilde\alpha$ are described as
a tuple of meromorphic functions on
$U \setminus\mathrm{Sing}\, D$.
Meromorphic (resp., holomorphic) functions can be extended
meromorphically (resp., holomorphically)
through an analytic subset of codimension at least two.
Hence $\alpha$ extends through $\mathrm{Sing}\, D$ to a meromorphic mapping
from $U$ to $\overline{A(M')}$.
{\it Q.E.D.}

We call a complex Lie group $T'$ a {\it quasi-torus} if
$T'$ carries an exact sequence of complex Lie groups similar to
(\ref{2.2}) with a complex torus $T$ instead of $A(M)$:
\begin{equation}
\label{2.4}
1 \to {(\mathbf{C}^*)^\rho}
 \to T' \stackrel{\pi}{\to} T \to 1.
\end{equation}
Then we can take a compactification $\overline {T'}$ of $T'$
as $\overline{A(M')}$ constructed above for $A(M')$, such that
$\bar\pi:\overline{T'} \to T$ is a holomorphic
$(\mathbf{P}^1)^\rho$-bundle
and the inclusion
$T'\hookrightarrow \overline{T'}$ is $T'$-equivariant.
Let $X$ be an irreducible complex subspace of
$\overline{T'}$ and set $X'=X\cap T'$.
Here we would like to discuss the stabilizer $\mathrm{St}(X)$ of $X$
defined as the subgroup consisting of those elements of $T'$
whose actions preserve the set $X$,
and the quotient $T'/\mathrm{St}(X)$.

Since $X'$ is dense in $X$, it is clear that
$$
\mathrm{St}(X)=\{ t \in T': t \cdot X' \subset X\},
$$
where ``$\cdot$'' stands for the action of $T'$ on $\overline{T'}$.
Using $X' \subset T'$, one may reformulate it as
$$
\mathrm{St}(X)=\{t \in T': t \cdot x \in X \hbox{ for }
\forall x \in X' \}
= T' \cap\ \left( \bigcap_{x \in X'} 
x^{-1}\cdot X \right),
$$
where $x^{-1}$ denotes the inverse of $x$ with
respect to the group law of $T'$.
The last description of $\mathrm{St}(X)$ implies 
that $\mathrm{St}(X)$ is Zariski-open in its closure in
$\overline{T'}$.
In particular it is a closed complex Lie subgroup of $T'$.
Moreover, $\mathrm{St}(X)$ is open and dense in the
compact analytic set $\overline{\mathrm{St}(X)}$
that is the closure of $\mathrm{St}(X)$ in $\overline{T'}$.
The proper mapping theorem implies that
$\bar\pi(\overline{\mathrm{St}(X)})$ is a compact analytic subset of $T$.
As a closure of the group
$\pi(\mathrm{St}(X))$, the compact set
$\bar\pi(\overline{\mathrm{St}(X)})$
must be a subgroup of $T$, too.
Now, $\pi(\mathrm{St}(X))$
is an open and dense subgroup
of the topological group $\bar\pi(\overline{\mathrm{St}(X)})$.
For any topological group, there never exists any open and
dense subgroup except the group itself.
Hence $\pi(\mathrm{St}(X))=\bar\pi(\overline{\mathrm{St}(X)})$
and
\begin{equation}
\label{2.5}
\pi(\mathrm{St}(X)) \hbox{ is a compact subgroup of } T.
\end{equation}

The intersection
$\overline{\mathrm{St}(X)} \cap (\overline{\mathbf{C}^*})^\rho$
is a closed analytic subset of 
$(\overline{\mathbf{C}^*})^\rho =(\mathbf{P}^1)^\rho$.
>From Chow's Theorem it follows that
$\overline{\mathrm{St}(X)} \cap (\overline{\mathbf{C}^*})^\rho$ 
is algebraic, and so is $\mathrm{St}(X)\cap (\mathbf{C}^*)^\rho$.
Thus $\mathrm{St}(X)\cap (\mathbf{C}^*)^\rho$
is an algebraic subgroup of $(\mathbf{C}^*)^\rho$.
Hence
\begin{equation}
\label{2.6}
(\mathbf{C}^*)^\rho/(\mathrm{St}(X)\cap (\mathbf{C}^*)^\rho)\cong
(\mathbf{C}^*)^\tau
\end{equation}
for some integer $\tau \geqq 0$.
This yields an short exact sequence for the quotient group
$T'/\mathrm{St}(X)$:
\begin{equation}
\label{2.7}
1 \to \mathrm{St}(X)\cap (\mathbf{C}^*)^\rho \cong (\mathbf{C}^*)^\tau  \to 
T'/\mathrm{St}(X) \to T/\pi(\mathrm{St}(X)) \to 1,
\end{equation}
where $T/\pi(\mathrm{St}(X))$ is a compact torus denoted by $T_1$.
We set $T'_1=T'/\mathrm{St}(X)$ which is again a quasi-torus.

Now, taking an algebraic equivariant compactification 
$(\mathbf{C}^*)^\tau \hookrightarrow (\overline{\mathbf{C}^*})^\tau=
(\mathbf{P}^1)^\tau$
of $(\mathbf{C}^*)^\tau$,
we obtain the $(\mathbf{P}^1)^\tau$-bundle over $T_1$
associated to the $(\mathbf{C}^*)^\tau$-bundle
$T'_1 \to T_1$, which yields a $T'_1$-equivariant compactification
$\overline{T'_1} $ of $T'_1$.

We now claim that the holomorphic projection $\lambda : T'\to T'_1$ extends
to a meromorphic mapping $\bar\lambda:\overline{T'} \to \overline{T'_1}$.
We consider the projection $\lambda:{T'} \to {T'_1}$ in two steps:
$$
T' \to T'/(\mathrm{St}(X)\cap (\mathbf{C}^*)^\rho) \to
T'/\mathrm{St}(X)=T'_1 .
$$
Setting $T'_2=T'/(\mathrm{St}(X)\cap (\mathbf{C}^*)^\rho)$,
we have a short exact sequence by (\ref{2.7})
\begin{equation}
\label{2.8}
1 \to (\mathbf{C}^*)^\tau \to T'_2 \to T \to 1.
\end{equation}
Using the above chosen algebraic compactification
$(\mathbf{C}^*)^\tau \hookrightarrow (\mathbf{P}^1)^\tau$,
we obtain a compactification $\overline{T'_2}$ of $T'_2$,
which is a $(\mathbf{P}^1)^\tau$-bundle over $T$
associated to the $(\mathbf{C}^*)^\tau$-principal bundle $T'_2 \to T$.
Both $T'$ and $T'_2$ are bundles over the same compact base manifold $T$.
Hence, in order to discuss the closure of the graph of the mapping
$T' \to T'_2$ in $\overline{T'}\times \overline{T'_2}$
we may restrict ourselves to the fibers over
the base manifold $T$.
But the fibers $(\mathbf{P}^1)^\rho$ and $(\mathbf{P}^1)^\tau$
are algebraic and the projection
$(\mathbf{C}^*)^\rho \to (\mathbf{C}^*)^\tau$
is an algebraic morphism of algebraic groups.
Hence the projection $(\mathbf{C}^*)^\rho \to (\mathbf{C}^*)^\tau$
extends to a meromorphic mapping
$(\mathbf{P}^1)^\rho \to (\mathbf{P}^1)^\tau$,
and consequently the projection $T' \to T'_2$ extends
to a meromorphic mapping
$\overline{T'} \to \overline{T'_2}$.

Next, let us discuss $T'_2 \to T'_1$.
Set
$$
C=\mathrm{St}(X)/((\mathbf{C}^*)^\rho \cap
\mathrm{St}(X))
$$
Then $C \cong \pi(\mathrm{St}(X))$ and
it is a compact torus by (\ref{2.5}).
Hence the projection $T'_2 \to T'_1$ is
the quotient of the manifold $T'_2$
by the action of the compact complex Lie group $C$.
In general, for actions of a {\it compact}
complex Lie group on a normal complex space,
there always exists a geometric quotient which
is again a normal complex space (see [H63], p.\ 358, Satz 20).
Since the projection ${T'_2} \to {T'_1}$
is $C$-equivariant, it extends to
a holomorphic mapping $\overline{T'_2} \to \overline{T'_1}$,
and moreover, the quotient $\overline{T'_2}/C$ is isomorphic
to $\overline{T'_1}$.

We have $X_1:=\bar\lambda(X)\subset \overline{T'_1}$ and
$X'_1:=\lambda(X')\subset T'_1$ such that
$X'_1=X_1 \cap T'_1$ and $\mathrm{St}(X'_1)=\{0\}$.
Then we deduce from the jet projection method used in
[N81], [\NO] and [N98] that $X'_1$ is of general type, and
specially, $X_1$ is a Moishezon space.
It follows from (\ref{2.7}) that
\begin{equation}
\label{2.9}
1 \to (\mathbf{C}^*)^\tau \to T'_1 \stackrel{\pi_{1_{}}}{\to}
 T_1 \to 1,
\end{equation}
where $\pi_1:T'_1 \to T_1$ extends to a holomorphic projection
$\bar\pi_1:\overline{T'_1} \to T_1$.
If $\pi_1(X_1)$ generates $T_1$, then $T_1$ is an Abelian variety,
and hence $T'_1$ is a so-called quasi-Abelian or semi-Abelian variety;
the sequence mappings are algebraic morphisms.

Summarizing the above obtained, we have
\begin{lemma}
\label{2.10}
Let $T'$ be a quasi-torus and let $X \subset \overline{T'}$ be
an irreducible complex subspace.
\begin{enumerate}
\item
The quotient space $T'/\mathrm{St}(X)$ is a quasi-torus,
and the quotient mapping
$\lambda:T' \to T'/\mathrm{St}(X)$ is holomorphic and
meromorphically extends to
$\bar\lambda:\overline{T'} \to \overline{T'/\mathrm{St}(X)}$.
\item
If $X'$ generates $T'$, then $T'/\mathrm{St}(X)$ is a quasi-Abelian
variety.
\end{enumerate}
\end{lemma}

\begin{theorem}
\label{2.11}
Let $M$ be a K\"ahler manifold and let $D$ be a hypersurface of $M$.
If $\bar q(M \setminus D)> dim M$, then the image
of an entire holomorphic curve $f:\mathbf{C} \to M\setminus D$
is contained in a proper analytic subset of $M$.
\end{theorem}

{\it Remark.}
For algebraic $M$ this was proved
by [N77], [N81], [\NO].

{\it Proof.}
Set $M'=M\setminus D$, and let $\alpha: M' \to A(M')$ be the
quasi-Albanese mapping.
It follows from Lemma \ref{2.3} that $\alpha$ extends meromorphically
to $\bar\alpha:M \to \overline{A(M')}$ between the compactifications.
Set $X=\bar\alpha(M)$ and $X'=X\cap A(M')$.
Then $X'$ generates $A(M')$ by virtue of the universal property
of the quasi-Albanese variety (Iitaka [I76]).
Set $A'_1=A(M')/\mathrm{St}(X)$.
It follows from the assumption that $\dim A'_1 > 0$.
Let $\lambda: A(M') \to A'_1$ be the quotient mapping
with the meromorphic extension
$\bar\lambda: \overline{A(M')} \to \overline{A'_1}$.
Set $X_1=\bar\lambda(X)$ and $X'_1=X_1 \cap A'_1$.
Then $\mathrm{St}(X_1)=\{0\}$ and $X'_1$ generates $A'_1$.
It follows from Lemma \ref{2.9} that $A'_1$ is quasi-Abelian
and $X'_1$ is an algebraic subvariety of $A'_1$.
Thus the case is reduced to algebraic one.
We know by [N81] that the composite
$\lambda\circ\alpha\circ f:\mathbf{C} \to X'_1$ is algebraically degenerate.
Thus the claim follows.
{\it Q.E.D.}

\section{Integral point sets}

In this section $K$ denotes a number field, and $M(K)$ denotes
a set of inequivalent non-trivial places (absolute values) of $K$
satisfying the product formula
$$
\prod_{v \in M(K)} | a |_v=1, \qquad a \in K .
$$
Let $M(K)_0$ denote the subset of all non-archimedean places of $M(K)$.
We have the ring $\mathcal{O}(K)$ of $K$-integers $a \in K$ defined by
$$
|a|_v \leqq 1, \quad \forall v \in M(K)_0.
$$
Let $S$ be a finite subset of $M(K)$ containing all infinite places.
We call $a \in K$ an $S$-{\it integer} if
$$
|a|_v \leqq 1, \qquad \forall v \in M(K)\setminus S.
$$
\begin{definition}{\rm
\label{3.0}
Let $R \subset M(K)_0$ be the complement of a finite subset,
possibly empty, of $M(K)_0$.
We call a collection $\{c_v\}_{v \in R}$ of real numbers
$c_v \in \mathbf{R}$ an $R$-{\it constant}
if but for finitely many $v \in R$,
$c_v=0$.}
\end{definition}

Let $V$ be a projective variety with a Cartier divisor $D$.
We are going to define the notion of a $(D,S)$-integral point set
of $V(K)$.
If $D$ is very ample, a set $A\subset V(K)$ shall be called
``$(D,S)$-integral point set'' if and only if there is
a proper embedding of $V\setminus D$ into an affine space
induced by sections of the line bundle associated to $D$
such that all the coordinates of all the points in the image of $A$
are $S$-integers.
To state a definition for arbitrary divisors
we need to introduce Weil functions (cf.\ [L83]).
These are actually non-archimedean fiber length functions (metrics)
of a line bundle.
We introduce them in a constructive way familiar with complex geometers.

Let $L=[D]$ be the line bundle determined by $D$.
Let $\{U_\lambda\}_{\lambda \in \Lambda}$ be an affine covering
of $V$ such that the restrictions $L|U_\lambda$ are trivial.
There is a system of transition functions $f_{\lambda\mu}$
on $U_\lambda \cap U_\mu$ associated with this trivialization of $L$.
Now first, we assume that $L$ has no base locus, so that
there is a basis $\{\sigma_i\}_{i=0}^N$ of $H^0(V,L)$ such that
$$
\Phi_L=(\sigma_0, \ldots , \sigma_N): V \to \mathbf{P}_K^N
$$
is a morphism.
Let $\sigma_i|U_\lambda$ be given by a regular function
$\sigma_{i\lambda} \in \mathcal{O}(U_\lambda)$.
Then
$$
\sigma_{i\lambda}=f_{\lambda\mu}\sigma_{i\mu}
\hbox{ on } U_\lambda \cap U_\mu .
$$
Let $\sigma=\{\sigma_\lambda\}\in H^0(V,L)$ be a section such that
the divisor $(\sigma)=D$, and define a $v$-length of
$\sigma$ by
$$
||\sigma (x)||_v
=\frac{|\sigma_\lambda(x)|_v}{\max_i |\sigma_{i\lambda}(x)|_v},
\qquad v \in M(K)_0, \quad x \in V(K),
$$
which is independent of $\lambda$.
Then we set
$$
n_v(\sigma)(x)=\log
\frac{1}{||\sigma(x)||_v},
\quad v \in M(K)_0,\; x \in V(K).
$$
Suppose that other choices of sections $\sigma_i$ and $\sigma$
similarly yield $n'(v,\sigma')$, where $(\sigma')=D$.
Then there is an $M(K)_0$-constant
$\{c_v\}$ such that
\begin{equation}
\label{3.1}
-c_v \leqq n_v(\sigma)(x) - n'_v(\sigma')(x) \leqq c_v.
\end{equation}
We define
\begin{equation}
\label{3.2}
n_v(D)(x)=n_v(\sigma)(x),\qquad
v \in M(K)_0,\; x \in V(K).
\end{equation}
Then $n_v(D)$ is well-defined up to a term bounded by an
$M(K)_0$-constant.

One can show that $n_v(D+D')=n_v(D)+n_v(D')$
for any two divisors $D$, $D'$ such that the associated
line bundles have no base locus.
Hence for an arbitrary divisor $D$ it is natural
to proceed as follows.
Let
$D_i\ (i=1,2)$ be effective divisors such that $D=D_1-D_2$ and
such that the line bundles $L_i$ associated to $D_i$ have both
no base locus.
Put
\begin{equation}
\label{3.3a}
n_v(D)=n_v(D_1)-n_v(D_2).
\end{equation}
This makes sense except where $n_v(D_1)=n_v(D_2)=-\infty$
(it happens in $D_1\cap D_2$).
However, equation (\ref{3.3a}) is equivalent to
a more explicit formula which makes sense even in $D_1\cap D_2$.
Namely,
taking a refinement of the above covering $\{U_\lambda\}$ if necessary,
we have a basis $\sigma_{jk}=\{\sigma_{jk\lambda}\}$
of $H^0(V, L_j)$ for $ j=1,2$,
such that $\sigma_{jk\lambda} \in \mathcal{O}(U_\lambda)$.
Then (\ref{3.3a}) yields
\begin{align}
\label{3.4}
n_v(D) = n_v( \sigma)=& \log
\frac{\max_k |\sigma_{1k\lambda}(x)|_v}{\max_l
|\sigma_{2l\lambda(x)}|_v}
\frac{1}{|\sigma_\lambda(x)|_v},\\
    \nonumber \\
& v \in M(K)_0,\; x \in V(K). \nonumber
\end{align}
for $x\not\in D_1\cap D_2$ and
we may take this as definition for $n_v(D)$ if $x\in D_1\cap D_2$.

Thus there is a well-defined (up to addition by a  term bounded by an
$M(K)_0$-constant) ``Weil function'' $n_v(D)$ for every
Cartier divisor.

\begin{definition}
\label{3.5}
{\rm
A subset $A \subset V(K)$ is called a $(D,S)$-{\it integral point set} if
there is an $M(K)\setminus S$-constant $\{c_v\}$ satisfying
$$
n_v(D)(x) \leqq c_v, \qquad \forall v \in M(K)\setminus S,\;
\forall x \in A.
$$
When $S$ contains no non-archimedian place,
$A$ is called simply a $D$-integral point set.
}
\end{definition}

Note the following:
\begin{enumerate}
\item
An arbitrary finite $A \subset V(K)$ out of $D$
is an $(D,S)$-integral point set.
\item
If $D$ is effective, then an arbitrary $(D,S)$-integral point set
is contained in $V(K)\setminus D$.
\item
The function $\{n_v(D)\}_{v \in M(K)\setminus S}$ is 
called a {\it Weil function}, and is an
analogue of the {\it counting function of poles} in the Nevanlinna theory
(cf.\ [V]).
This is the reason why we use the notation $n_v(D)$ here.
\item
For two divisors $D, D'$ we have a natural addition formula
$$
n_v(D+D')=n_v(D)+n_v(D').
$$
\end{enumerate}

For example, we take a very ample divisor $D$ of $V$.
Let $\{\sigma_0, \sigma_1, \ldots , \sigma_N\}$ be a basis
of $H^0(V, [D])$ such that $(\sigma_0)=D$.
Setting $\phi_i=\sigma_i/\sigma_0, 1 \leqq i \leqq N$,
we have an affine imbedding
$$
\phi=(\phi_1, \ldots,\phi_N): V \setminus D \hookrightarrow \mathbf{A}_K^N.
$$
By definition
$$
n_v( D)(x)=\log \max\{ |\phi_i(x)|_v \}_{i=1}^N, \quad x \in
V(K)\setminus D.
$$
Therefore a subset $A \subset V(K)\setminus D$ is a
$(D,S)$-integral point set if and only if
there is an integer $c \in \mathcal{O}(K)$
such that all $c\phi_i(x)$ are $S$-integers for
$1 \leqq i \leqq N$ and $x \in A$.

\section{Proofs}

We will simultaneously give the proofs of our statements
for holomorphic curves and integral point sets, where
the arguments are analogous.

{\it (a) Proof of Theorem \ref{1.1}, (i).}

Let $W$ and $D_i, 1 \leqq i \leqq l$, be as in Theorem \ref{1.1}.
Set $l_0=\#\{W \cap D_i\not=W\}$.
Changing the indices, we have
$$
\{W \cap D_i\}=\{W \cap D_1,\ldots,W\cap D_{l_0}\},
$$
where $W\cap D_i, 1 \leqq i \leqq l_0$, are distinct hypersurfaces
of $W$.
Let $\gamma:\tilde W \to W$ be a desingularization, and set
$$
\tilde D_i=\gamma^{-1}(W\cap D_i), \quad 1 \leqq i \leqq l_0, \qquad
\tilde D=\sum_{i=1}^{l_0} \tilde D_i.
$$
Then $\tilde D_i, 1 \leqq i \leqq l_0$, are distinct.
Set
$$
r=\rankcDi .
$$
Since the Chern classes
$c_1(D_i), 1\leqq i \leqq l$, generate a subgroup of $H^2(M,\mathbf{R})$
of $\mathbf{Z}$-rank $r$, there are at least $l_0-r$ independent
$\mathbf{Z}$-linear relations between $c_1(D_i), 1 \leqq i \leqq l_0$.
Let $\sum_{i=1}^{l_0} a_{ki} c_1(D_i)=0, 1 \leqq k \leqq l_0-r$,
$a_{ki} \in \mathbf{Z}$, be such relations.
By changing the indices and by elementary transformations of
matrices one may assume that
\begin{alignat}{2}
\label{4.1}
a_{ii} &\not= 0, & &1 \leqq i \leqq l_0-r, \\
a_{ki} &=0, &\qquad &1 \leqq k\not= i \leqq  l_0-r.
\nonumber
\end{alignat}
Since $M$ is K\"ahler, there are multiplicative meromorphic functions
$\vartheta_k, 1 \leqq k \leqq l_0-r$, such that
the residue of logarithmic 1-forms $d\log \vartheta_k$ is
$\sum_{i=1}^{l_0}a_{ki}D_i$ (cf., e.g., [N77], [N81]).
Since $D_i\cap W, 1 \leqq i \leqq l_0-r$, are distinct,
logarithmic 1-forms $\gamma^* d\log \vartheta_k, 1 \leqq k \leqq l_0-r$,
are linearly independent over $\mathbf{C}$.
Thus we have
\begin{equation}
\label{4.2}
\bar q(\tilde W \setminus \tilde D)\geqq q(\tilde W)+ l_0-r=q(W)+ l_0- r.
\end{equation}
Since $f:\mathbf{C} \to W\setminus D$ lifts to
a holomorphic curve $\tilde f:\mathbf{C} \to \tilde W\setminus \tilde D$
with Zariski dense image, Theorem \ref{2.11} and (\ref{4.2}) imply that
\begin{equation}
\label{4.3}
q(W)+ l_0-r \leqq \dim W.
\end{equation}
This proves (i) of Theorem \ref{1.1}.

{\it (b) Proofs of Theorem \ref{1.1}, (ii), and Theorem \ref{1.5}.}

We first prepare a result from the Diophantine approximation.
Let $K$, $M(K)$, $S$, and $V$ be as in \S3.

\begin{lemma}
\label{4.3.1}
Let $\{D_i\}_{i=1}^l$ be a family of hypersurfaces of $V$ in general
position.
Let $W \subset V$ be a subvariety such that there is a
Zariski dense 
$(\sum_{D_i \not\supset W} D_i\cap W,S)$-integral
point set $A$ in $W$.
Then we have that
$$
\#\{W \cap D_i\not=W\} + q(W) \leqq \dim W + \rankcDi.
$$
\end{lemma}

This is almost the same as Vojta [V96], Corollary 0.3,
and follows from [V96], Theorem 0.2 and (\ref{4.2}).

To prove Theorem \ref{1.1}, (ii)
and Theorem \ref{1.5}, we count the number $l_0=\#\{W \cap D_i\not=W\}$.
In the sequel, we will take a finite extension of $K$ as much as we
need, but this does not weaken our claim.

Set $n=\dim W$.
For every $1\leqq j \leqq l_0$, let $s_j$ denote the number
of hypersurfaces $D_i$ with $D_i\cap W=D_j\cap W$.
Furthermore, let $s_0$ denote the number of divisors $D_i$
with $W\subset D_i$.
Changing the indices, we may assume that
\begin{equation}
\label{4.4}
s_1\geqq s_2 \geqq \cdots \geqq s_{l_0}.
\end{equation}
Let $l_1=\min\{n,l_0\}$.
Since the divisors are ample and in general position, it is clear
that the intersection 
$W \cap (\bigcap_{i=1}^{l_1} D_i)$ is non-empty.
Since the intersection of more than $m$ of the divisors $\{D_i\}_1^l$
must be empty, it follows that
\begin{equation}
\label{4.5}
\sum_{k=0}^{l_1} s_k \leqq m
\end{equation}
For $l_0 \leqq n$ this yields $l \leqq m$.

For $l_0 > n$ we have $l_1=n$.
By (\ref{4.3}) we have
$$
r-q(W)\geqq l_0 -n >0.
$$
It follows from  (\ref{4.4}) that
$$
\frac{1}{l_0}\sum_{k=1}^{l_0} s_k \leqq \frac{1}{n} \sum_{k=1}^{n} s_k,
$$
so that
\begin{equation}
\label{4.6}
\sum_{k=1}^{l_0} s_k \leqq \frac{l_0}{n} \sum_{k=1}^{n} s_k.
\end{equation}
Since $s_0 \leqq \frac{l_0}{n}s_0$, it follows from 
(\ref{4.6}), (\ref{4.5}) and (\ref{4.3}) that
$$
l=\sum_{k=0}^{l_0} s_k \leqq \frac{l_0}{n} \sum_{k=0}^n s_k \leqq
\frac{l_0 m}{n} \leqq m \left(1 + \frac{r-q(W)}{n}\right).
$$

{\it (c) Proof of Corollaries \ref{1.4} and \ref{1.6}.}

Due to M. Green [G77], there is a stratified version
of Brody's theorem (cf.\ [Ko98], Theorem (3.6.13) for this generalized
form):
\begin{lemma}
\label{4.7}
Let $X$ be a compact complex space and let $\{E_i\}_{i\in I}$ be
a family of Cartier hypersurfaces of $X$. Assume that
for every subset $\emptyset\subseteq J\subseteq I$,
every holomorphic curve
$$
f:\mathbf{C} \to
\bigcap_{j\in J} E_j \setminus \bigcup_{i\in I\setminus J} E_i
$$
is reduced to a constant mapping, where $\bigcap_{j\in\emptyset}E_j=X$.
Then $X \setminus\bigcup_i E_i$ is complete hyperbolic
and hyperbolically imbedded into $X$.
\end{lemma}

Suppose that (i) of Corollary \ref{1.4} or \ref{1.6} is false.
In case of Corollary \ref{1.4},
it follows from Lemma \ref{4.7} that for some subset
$\emptyset\subseteq J\subseteq I=\{1, \ldots, l\}$,
there exists a non-constant holomorphic curve
$$
f:\mathbf{C} \to
\bigcap_{i\in J} D_i \setminus \bigcup_{i\in I\setminus J} D_i
$$
Let $W$ be the Zariski closure of $f(\mathbf{C})$ in $M$.
Or, in case of Corollary \ref{1.6}, there exists
an infinite $(\sum_{i=1}^l D_i,S)$-integral point set $A$ of $V(K)$.
Let $W$ be a positive dimensional irreducible component
of the Zariski closure of $A$ in $V$.
Then Theorem \ref{1.1}, (ii), or Theorem \ref{1.5} implies that
$l \leqq m(1+\mathrm{rank}_{\mathbf{Z}}\, \mathrm{NS}(M\hbox{ or }V))$;
this is a contradiction.

To prove (ii) and (iii) of Corollaries \ref{1.4} and \ref{1.6}
we first note that for $\mathbf{P}^m(\mathbf{C})$
we get
\begin{equation}
\label{4.8}
\rankcDi=1.
\end{equation}
Then we have (ii) by arguments similar to the proofs of
Theorem \ref{1.1}, (ii)/Theorem \ref{1.5} and
(i) of Corollary \ref{1.4}/\ref{1.6}.

Claim (iii) is a direct consequence of Theorem \ref{1.1}, (ii)
or of Theorem \ref{1.5} and (\ref{4.8}).

{\it (d) Proof of Proposition \ref{1.4}.}

First, we note that in Lemma \ref{4.7} $f$ may be assumed
to be of order at most 2;
hence, in the present case $f$ is a translate of a
one-parameter subgroup of $M$.
Furthermore, by [A72] such $f$, if exists, must take values in
any ample divisor of $M$.

Assume that (i) is false.
Then there is a translate $W$ of a positive dimensional closed
subgroup $B$ of $M$ such that $W \subset \bigcap_{i=1}^l D_i$.
Thus, $m_0 \leqq \dim W \leqq \dim M -l$,
which contradicts the assumption.

In the case of (ii), we have a subset $J$ of the
index set $I=\{1,\ldots, l\}$ such that
$$
W \subset \bigcap_{i \in J} D_i, \qquad
W \cap \bigcup_{i\in I\setminus J} D_i=\emptyset.
$$
The first inclusion relation implies
\begin{equation}
\label{4.9}
m_0 \leqq \dim B=\dim W \leqq \dim M - \# J.
\end{equation}
By making use of the quotient mapping $\pi:M \to M/B$ one sees that
$D_i$ with $i \not\in J$ are $B$-invariant.
Since $D_i, i \in I$, are in general position,
\begin{equation}
\label{4.10}
m_0 \leqq \dim B \leqq \dim M - \# I \setminus J.
\end{equation}
It follows from (\ref{4.9}) and (\ref{4.10}) that
$2m_0 \leqq 2 \dim M - l$; this contradicts the assumption.

\begin{remark}
\label{4.12}{\rm
(i) The inequalities obtained in Theorems \ref{1.1}, \ref{1.5},
Corollaries \ref{1.3} and \ref{1.6} are optimal in general.
To confirm it it suffices to show that the estimate in (ii) of
Corollaries \ref{1.3} and \ref{1.6},
$$
\dim W \leqq \frac{m}{l-m},
$$
where $M=V=\mathbf{P}^m$, is optimal.
Take an arbitrary integer $l$ such that $m < l \leqq 2m$.
Let $m_0$ be the largest integer which does not exceed $m/(l-m)$.
For an arbitrarily given family $\{D_i\}_{i=1}^l$
of hyperplanes of $\mathbf{P}^m_{\mathbf{Q}}$
in general position,
Fujimoto [Fu72] and Green [G72] constructed
an imbedding of the $m_0$-dimensional multiplicative torus
$\mathbf{G}_{\mathbf{Q}}^{m_0}$ of non-zero
rational numbers defined over $\mathbf{Q}$,
$$
\mathbf{G}_{\mathbf{Q}}^{m_0} \hookrightarrow
\mathbf{P}^m_{\mathbf{Q}}\setminus \bigcup_{i=1}^l D_i.
$$
By taking a finite extension $K$ of $\mathbf{Q}$ or
$\mathbf{C}$, it is just easy to have such 
a variety $W$ of dimension $m_0$ as it contains
a Zariski dense $(\sum D_i\cap W, S)$-integral point set of $W(K)$,
or over $\mathbf{C}$ it admits an
entire holomorphic curve with Zariski dense image.

(ii) In [MN96] we constructed a number of irreducible
hyperbolic hypersurfaces $D$ of $\mathbf{P}^m_{\mathbf{C}}$
whose complements are hyperbolic and hyperbolically
imbedded into $\mathbf{P}^m_{\mathbf{C}}$.
One may define those $D$ over any number fields.
It is an interesting problem to study
$(D,S)$-integral point set of the complement of $D$.
In fact, we proved in [N97] the corresponding finiteness
theorems for function fields and $S$-unit points over
a number field.}
\end{remark}

\section*{\centerline{\normalsize\it References}}

{\small\parindent35pt
\begin{itemize}
\setlength{\itemsep}{-3pt}
\item[{[Ax72]}]
J. Ax,
Some topics in differential algebraic geometry II,
Amer.\ J. Math.\
{\bf 94}
(1972),
1205-1213.

\item[{[B84]}]
V.A. Babets,
Picard-type theorems for holomorphic mappings,
Siberian Math.\ J.
{\bf 25}
(1984),
195-200.
\item[{[Bl26]}]
A. Bloch,
Sur les syst\`emes de fonctions uniformes satisfaisant \`a l'\'equation d'une vari\'et\'e alg\'ebrique dont l'irr\'egularit\'e d\'epasse la dimension,
J. Math.\ Pures Appl.\
{\bf 5} (1926), 9-66.
\item[{[D71]}]
P. Deligne,
Th\'eroie de Hodge, II,
I.H.E.S. Publ.\ Math.\ {\bf 40} (1971), 5-57.
\item[{[F91]}]
G. Faltings,
Diophantine approximation on abelian varieties,
An..\ Math.\ {\bf 133} (1991), 549-576.
\item[{[Fu72]}]
H. Fujimoto,
Extension of the big Picard's theorem,
Tohoku Math.\ J.
{\bf 24}
1972,
415-422.
\item[{[G72]}]
M. Green,
Holomorphic maps into complex projective space,
Trans.\ Amer.\ Math.\ Soc.\
{\bf 169}
(1972),
89-103.
\item[{[G77]}]
M. Green,
The hyperbolicity of the complement of $2n + 1$ hyperplanes in general position in $\mathbf{P}_n$, and related results,
Proc.\ Amer.\ Math.\ Soc.\
{\bf 66}
(1977),
109-113.
\item[{[G78]}]
M. Green,
Holomorphic maps to complex tori,
Amer.\ J. Math.\
{\bf 100}
(1978),
615-620.
\item[{[H63]}]
H. Holmann,
Komplexe R\"aume mit komplexen Transformationsguppen,
Math.\ Ann.\ {\bf 150} (1963), 327-360.
\item[{[I76]}]
S. Iitaka,
Logarithmic forms of algebraic varieties,
J. Fac.\ Sci.\ Univ.\ Tokyo, Sect.\ IA {\bf 23} (1976), 525-544.
\item[{[K80]}]
Y. Kawamata,
On Bloch's conjecture,
Invent. Math. {\bf 57} (1980), 97-100.
\item[{[Ko70]}]
S. Kobayashi, Hyperbolic Manifolds and Holomorphic Mappings, 
Marcel Dekker, New York, 1970.
\item[{[Ko98]}]
S. Kobayashi, Hyperbolic Complex Spaces, Grundlehren der mathematischen
Wissenschaften {\bf 318}, Springer-Verlag, Berlin-Heidelberg, 1998.
\item[{[L83]}]
S. Lang,
Fundamentals of Diophantine Geometry,
Springer-Verlag, New York-Berlin-Heidelberg-Tokyo, 1983
\item[{[L87]}]
S. Lang,
Introduction to Complex Hyperbolic Spaces,
Springer-Verlag, New York-Berlin-Heidelberg, 1987.
\item[{[L91]}]
S. Lang,
Number Theory III, Encycl.\ Math.\ Sci.\ vol.\ {\bf 60},
Springer-Verlag,
Berlin-Heidelberg-New York-London-Paris-Tokyo-Hong Kong-Barcelona, 1991.
\item[{[MN96]}]
K. Masuda and J. Noguchi,
A construction of hyperbolic hypersurfaces of
$\mathbf{P}^n (\mathbf{C})$,
Math.\ Ann.\ {\bf 304} (1996), 339-362.
\item[{[N77]}]
J. Noguchi,
Holomorphic curves in algebraic varieties,
Hiroshima Math.\ J.
{\bf 7}
(1977),
833-853.
\item[{[N80]}]
J. Noguchi,
Supplement to ``Holomorphic curves in algebraic varieties",
Hiroshima Math.\ J.
{\bf 10}
(1980),
229-231.
\item[{[N81]}]
J. Noguchi,
Lemma on logarithmic derivatives and holomorphic curves in algebraic varieties,
Nagoya Math.\ J.
{\bf 83}
(1981),
213-233.
\item[{[N91]}]
J. Noguchi,
Hyperbolic Manifolds and Diophantine Geometry,
Sugaku Exposition Vol.\
{\bf 4}
pp.\ 63-81,
Amer.\ Math.\ Soc., Rhode Island,
1991.
\item[{[N92]}]
J. Noguchi,
Meromorphic mappings into compact hyperbolic complex spaces and geometric
Diophantine problems,
International J. Math.\
{\bf 3}
(1992),
277-289.
\item[{[N95]}]
J. Noguchi,
A short analytic proof of closedness of logarithmic forms,
Kodai Math.\ J. {\bf 18} (1995), 295--299.
\item[{[N96]}]
J. Noguchi,
On Nevanlinna's second main theorem,
Geometric Complex Analysis,
Proc.\ the Third International Research Institute, Math.\ Soc.\ Japan,
Hayama, 1995, pp.\ 489-503, World Scientific, Singapore, 1996.
\item[{[N97]}]
J. Noguchi,
Nevanlinna-Cartan theory over function fields and a Diophantine equation,
J. reine angew.\ Math.\ {\bf 487} (1997), 61-83.
\item[{[N98]}]
J. Noguchi,
On holomorphic curves in semi-Abelian varieties,
Math.\ Z. {\bf 228} (1998), 713-721.
\item[{[\NO]}]
J. Noguchi and T. Ochiai,
Geometric Function Theory in Several Complex Variables,
Japanese edition, Iwanami, Tokyo, 1984;
English Translation, Transl.\ Math.\ Mono.\ {\bf 80},
Amer.\ Math.\ Soc., Providence, Rhode Island,
1990.
\item[{[O77]}]
T. Ochiai,
On holomorphic curves in algebraic varieties with ample irregularity,
Invent.\ Math.\ {\bf 43} (1977), 83-96.
\item[{[RW91]}]
M. Ru and P.-M. Wong,
Integral points of $\mathbf{P}^n-\{2n+1\hbox{ hyperplanes in general position}\}$,
Invent.\ Math.\ {\bf 106} (1991), 195-216.
\item[{[S-Y96]}]
Y.-T. Siu and S.-K. Yeung,
A generalized Bloch's theorem and the hyperbolicity of the complement
of an ample divisor in an Abelian variety,
Math.\ Ann.\ {\bf 306} (1996), 743-758.
\item[{[V96]}]
P. Vojta,
Integral points on subvarieties of semiabelian varieties, I,
Invent.\ Math.\ {\bf 126} (1996), 133-181.

\end{itemize}
}
\bigskip

\rightline{Graduate School of Mathematical Sciences}
\rightline{University of Tokyo}
\rightline{Komaba, Meguro,Tokyo 153-8914}
\rightline{e-mail: noguchi@ms.u-tokyo.ac.jp} 
\medskip
\rightline{Mathematisches Institut}
\rightline{Rheinsprung 21}
\rightline{CH--4053 Basel}
\rightline{Switzerland}
\rightline{e-mail: jwinkel@member.ams.org} 
\rightline{URL: http://www.cplx.ruhr-uni-bochum.de/$\tilde{\ }$jw/index-e.html}

\end{document}